# Close Approximations for Daublets and their Spectra

V. Vermehren V. and J. E. Wesen
State University of Amazon
Electrical Engineering Department
Manaus – Amazon, Brazil
vvv@netium.com.br and ednelson.wesen@gmail.com

H. M. de Oliveira.
Federal University of Pernambuco
Signal Processing Group
Recife – Pernambuco, Brazil
hmo@ufpe.br

*Abstract*—This paper offers a new regard on compactly supported wavelets derived from FIR filters. Although being continuous wavelets, analytical formulation are lacking for such wavelets. Close approximations for daublets (Daubechies wavelets) and their spectra are introduced here. The frequency detection properties of daublets are investigated through scalograms derived from these new analytical expressions. These near-daublets have been implemented on the Matlab$^{TM}$ wavelet toolbox and a few scalograms presented. This approach can be valuable for wavelet synthesis from hardware or for application involving continuous wavelet-based systems, such as wavelet-OFDM.

## I. INTRODUCTION

Wavelets have lately gained prolific applications throughout an amazing number of areas, especially in Physics and Engineering [1]. Both Continuous and Discrete Wavelet transforms (CWT and DWT, respectively) have emerged as a definitive tool of signal processing analysis and have proven to be superior to classical Fourier analysis in many situations [1–3]. Continuous wavelet transform furnishes precise information on the local and global regularity [4] is very robust to spurious contamination of the signal. It has also been used in damping identification [5], analysis of seismic signals [6] and detection of singularities [4], among many other applications. A number of continuous wavelets do have analytical close expression, such as Morlet, Meyer, Mexican hat, Shannon, beta etc. [3, 7]. However, an astonishing, but well established feature of wavelets is that in many cases the signal analysis can be carried out without knowing the analytical expression of the mother wavelet. This fact is particularly true for most compactly supported wavelets [7]. The plot or the generation of the waveform of both the wavelet and the scaling function is always performed by an iterative numeric approach (cascade algorithm). How could one generate a db waveform without exploiting numeric successive approximations? This paper offers an (approximate) analytical solution, which turns straightforward to handle with db-type wavelets (e.g. coiflets). The paper is organized as follows. Section II introduces the idea of inharmonic series as some sort of generalization of standard Fourier series. The derivation of analytical expressions for *db*N wavelets is carried out in Section III, and tight continuous approximations are supplied both in time and frequency domain. The investigation of the frequency detection properties of *db*N wavelets were addressed in Section IV through scalograms. Although scalograms suggest a poor detection feature, it is shown that *db*N actually have good frequency detection properties. Section V presents some potential applications of the analytical formulation of Daubechies wavelets in the framework of Telecommunication. Finally, some concluding remarks are presented in Section VI.

## II. INHARMONIC SERIES ANALYSIS

Let us consider a compactly supported wavelet $\psi(t)$ defined in the interval $[0, T]$, where

$$T := \text{length}(\text{supp}\,\psi(t)) \tag{1}$$

is the length of the support of the wavelet [1]. The classical model (Fourier) to represent a signal [2] within a window of length $T$ is given by:

$$\psi(t) = \sum_{k=1}^{\infty} A_k \sin\left(\tfrac{2\pi}{T} k t + \theta_k\right). \tag{2}$$

There are harmonic components with a fundamental $2\pi/T$. Instead of this standard well-known model, we deal with an unbounded supported Fourier-like wavelet representation,

$$\psi_{long}(t) = \sum_{k=1}^{\infty} A_k \sin\left(\omega_k t + \theta_k\right), \tag{3}$$

where the frequencies $\omega_k$ are no more harmonic frequencies. Within the wavelet support, $\psi_{long}(t) = \psi(t)$, where $\psi(t)$ is the wavelet. Imposing now the oscillation condition to the wavelet [7], it follows that

$$\int_0^T \psi(t)dt = \int_0^T \psi_{long}(t)dt, \tag{4}$$

which yields

$$\sum_{k=1}^{\infty} A_k \int_0^T \sin(\omega_k t + \theta_k)dt = 0. \tag{5}$$

In order to assure a zero-mean for the wavelet, a possible condition is that each integral in the former equation vanishing, giving

$$\int_0^T \sin(\omega_k t + \theta_k)dt = \int_{\theta_k}^{\omega_k T + \theta_k} \sin(\zeta)d\zeta = 0. \tag{6}$$

This can be rewritten as

$(1 - \cos(\omega_k T_k)).(\cos(\theta_k)) + \sin(\omega_k T_k).(\sin(\theta_k)) = 0, \forall k$.

Here, some phase constraints are introduced, namely

$$\theta_k \neq \tfrac{2\pi}{T_k} l \quad \text{and} \quad l \in Z.$$

This work was partially supported by the Brazilian National Council for Scientific and Technological Development (CNPq).

From that above mentioned condition, it can be derived the following relationships:

1) Finding out the phase angles $\theta_k$ in terms of the frequencies components $\omega_k$,

$$\theta_k = -tg^{-1}\left(\frac{1-\cos(\omega_k T_k)}{\sqrt{1-\cos^2(\omega_k T_k)}}\right). \quad (7)$$

2) Another very interesting view is to determine the frequencies $\omega_k$ in terms of the $\theta_k$. Hence it follow that:

$$(1+tg^2(\theta_k))x_k^2 - 2x_k + (1-tg^2(\theta_k)) = 0, \quad (8)$$

where $x_k := \cos(\omega_k T_k)$.

The (non trivial) solution of interest is

$$\omega_k = \pm\frac{1}{T_k}\cos^{-1}\left(\frac{1-tg^2(\theta_k)}{1+tg^2(\theta_k)}\right). \quad (9)$$

In the sequel, a new analytical model is proposed to represent compactly supported wavelets in the interval [0, T]. This will be referred to as inharmonic series of a wavelet,

$$\psi_{long}(t) = \sum_{k=1}^{\infty} A_k \sin\left[\left(\frac{2\pi}{T}k \pm \frac{1}{T_k}\cos^{-1}\left(\frac{1-tg^2(\theta_k)}{1+tg^2(\theta_k)}\right)\right)t + \theta_k\right]. \quad (10)$$

The non-equally spaced frequencies of the decomposition are

$$\frac{2\pi}{T}k \pm \frac{1}{T_k}\cos^{-1}\left(\frac{1-tg^2(\theta_k)}{1+tg^2(\theta_k)}\right), \quad (11)$$

which are rather similar to the harmonic perturbations mentioned in [8] and derived from the wavelet oscillatory condition:

$$\int_0^T \psi(t)dt = 0. \quad (12)$$

The frequencies present in this series expansion are no more harmonic frequencies such as ($\omega_0, 2\omega_0, 3\omega_0, \ldots, \omega_0 := 2\pi/T$), but rather have perturbations inserted so as to guarantee a zero mean for the wavelet. This can be called an inharmonic decomposition of a wavelet.

III. THE ANALYTICAL (CLOSE) EXPRESSIONS FOR TIME AND SPECTRUM DOMAIN OF DAUBECHIES WAVELETS

In this section we derive approximations for both the dbN wavelets and their spectra.

A. Time Domain Analysis

The analytical (close) expression to implement a fairly accurate daublet in the time domain is

$$\psi_{long\,db6}(t) = \sum_{k=1}^{\infty} a_k \sin(b_k t + c_k). \quad (13)$$

This expression is a sum of *sin* functions with coefficients of amplitude, frequency and phase, assembled through nonlinear least square method with Levenberg-Marquardt robust algorithm [9, 10]. Typically just eight to ten terms are enough to get the *Daubechies* (*db*N) wavelets with 99.99% of confidence level.

Table I shows the coefficients $a_k$, $b_k$, $c_k$ required to perform close approximations for *db*4, *db*6 and *db*8 daublets. Description of such wavelets can be found in [1, 7].

Identical analytical (close) formulation can be used to perform approximations for daublets scale function ($\phi(t)$), by just using another set of coefficients.

The coefficients associated with Daubechies scaling function are showed in Table II. It can be remarked that they maintain the same features as previously reported to generate approximations for compactly supported wavelet functions.

TABLE I. COMPUTED COEFFICIENT TO APROXIMATE DAUBECHIES WAVELET FUNCTIONS ACCORDING TO (1) FOR $\psi_{DB4}$, $\psi_{DB6}$, AND $\psi_{DB8}$ WAVELETS.

| | $\psi_{db4}$ | | |
|---|---|---|---|
| k | $a_k$ | $b_k$ | $c_k$ |
| 1 | 0.3452 | 4.586 | -2.316 |
| 2 | 0.2783 | 3.460 | 1.413 |
| 3 | 0.3015 | 5.770 | -0.373 |
| 4 | 0.2129 | 6.960 | -4.943 |
| 5 | 0.1293 | 2.414 | -1.794 |
| 6 | 0.1120 | 8.161 | -3.225 |
| 7 | 0.0295 | 9.366 | -7.567 |
| 8 | 0.0223 | 1.372 | 1.102 |

| | $\psi_{db6}$ | | |
|---|---|---|---|
| k | $a_k$ | $b_k$ | $c_k$ |
| 1 | 0.2623 | 4.850 | -1.655 |
| 2 | 0.2520 | 3.993 | 3.014 |
| 3 | 0.2287 | 5.724 | 0.649 |
| 4 | 0.1778 | 3.197 | 0.649 |
| 5 | 0.1729 | 6.590 | -5.635 |
| 6 | 0.1098 | 7.459 | -4.613 |
| 7 | 0.0820 | 2.436 | 4.117 |
| 8 | 0.0504 | 8.333 | -9.828 |

| | $\psi_{db8}$ | | |
|---|---|---|---|
| k | $a_k$ | $b_k$ | $c_k$ |
| 1 | -0.2054 | 5.066 | -17.48 |
| 2 | 0.1334 | 3.116 | -6.671 |
| 3 | 0.1926 | 3.720 | -10.66 |
| 4 | -0.0622 | 2.532 | -12.39 |
| 5 | 0.2145 | 4.379 | -15.39 |
| 6 | 0.1768 | 5.750 | -26.04 |
| 7 | 0.1360 | 6.419 | -25.21 |
| 8 | -0.0917 | 7.081 | -27.53 |
| 9 | -0.0468 | 7.740 | -32.94 |

TABLE II. COMPUTED COEFFICIENTS TO APPROXIMATE DAUBECHIES SCALES FUNCTIONS ACCORDING TO (10) FOR $\phi_{DB4}$, $\phi_{DB6}$ AND $\phi_{DB8}$ WAVELETS.

| | $\phi_{db4}$ | | |
|---|---|---|---|
| k | $a_k$ | $b_k$ | $c_k$ |
| 1 | 0.3762 | 0.672 | 0.171 |
| 2 | 0.2113 | 3.226 | -2.404 |
| 3 | 0.3900 | 1.204 | 0.939 |
| 4 | 0.0770 | 4.193 | 2.098 |
| 5 | 0.2661 | 2.384 | -1.379 |
| 6 | 0.0081 | 5.586 | -1.379 |
| 7 | 0.0226 | 8.537 | -1.184 |
| 8 | 0.0205 | 9.424 | 3.346 |

| | $\phi_{db6}$ | | |
|---|---|---|---|
| k | $a_k$ | $b_k$ | $c_k$ |
| 1 | 0.2247 | 0.648 | 1.540 |
| 2 | 0.1244 | 1.323 | -0.241 |
| 3 | 0.3148 | 2.333 | -1.329 |
| 4 | 0.0111 | 0.032 | 0.8670 |
| 5 | 0.3007 | 2.084 | -2.628 |
| 6 | 0.0489 | 4.087 | -5.627 |
| 7 | 0.1224 | 3.019 | 2.881 |
| 8 | 0.0935 | 3.728 | 0.425 |
| 9 | 0.0296 | 0.338 | 0.208 |
| 10 | 0.2088 | 0.342 | 0.735 |

| | $\phi_{db8}$ | | |
|---|---|---|---|
| k | $a_k$ | $b_k$ | $c_k$ |
| 1 | 0.1417 | -0.004 | 1.617 |
| 2 | 0.1214 | 1.697 | -1.584 |
| 3 | 0.1480 | 2.174 | -2.969 |
| 4 | 0.1840 | -0.271 | 0.929 |
| 5 | 0.1603 | 2.544 | -3.420 |
| 6 | 0.1057 | 2.934 | -4.170 |
| 7 | 0.1136 | 3.586 | -1.468 |
| 8 | 0.0877 | 3.759 | -0.154 |
| 9 | 0.1234 | 0.907 | -0.523 |
| 10 | 0.1419 | 1.239 | -0.457 |

Now, taking $\psi_{db4}$ coefficients as a reference and setting $b_k$ in ascendant order, we can apply (10) so as to retrieve the corrections terms around $k\omega_0$ corresponding to $db4$ inharmonics (9). Table III displays these correction terms of fluctuation on the harmonic components of the series. The values of $k$, in contrast with those of Table I, now correspond to the order of the harmonic components in the series, with $w_0 = 2\pi/7$.

### B. Spectrum Domain Analysis.

The analytical (close approximations) expression in the spectrum domain for Daubechies wavelets is given by:

$$\Psi_{long}(\omega) = j\pi \sum_{k=1}^{\infty} a_k \left[ e^{-jc_k} \delta(\omega + b_k) - e^{jc_k} \delta(\omega - b_k) \right]. \quad (14)$$

This expression is the Fourier pair with the *inharmonic series* $\psi_{long}(t)$. Clearly, merely Dirac impulses appear representing *sin* functions on the spectrum of the periodic approximation of the analyzed wavelet [11]. As an example, let us consider the inharmonic expansion of a standard wavelet such as $db6$. The spectrum of the unbounded periodic approximated wavelet long $db6$ is shown in Fig.1.

TABLE III. COMPUTED $\omega_k$ COEFFICIENTS REQUIRED TO APPROXIMATE DAUBECHIES WAVELET ACCORDING TO (10) FOR THE $\psi_{DB4}$ WAVELET.

| k | $\psi_{db4}$ $b_n$ | harm. $k.w_0$ | inharm. $\omega_k$ |
|---|---|---|---|
| 1 | 1.372 | 0.898 | 0.474 |
| 2 | 2.414 | 1.795 | 0.619 |
| 4 | 3.460 | 3.590 | -0.130 |
| 5 | 4.586 | 4.488 | 0.096 |
| 6 | 5.770 | 5.386 | 0.384 |
| 7 | 6.960 | 6.283 | 0.677 |
| 9 | 8.161 | 8.078 | 0.083 |
| 10 | 9.366 | 8.976 | 0.390 |

In order to determine the "true" $db6$ spectrum, the wavelet $\psi_{longdb6}(t)$ must be confined in its support, according to:

$$\psi_{db6}(t) = \psi_{long\,db6}(t).\Pi\left(\frac{t-5.5}{11}\right), \quad (15)$$

where $\Pi$ denotes the standard gate function and T:=2N-1=11 is the length of the support of the db6 wavelet [1].

The spectrum of $\psi_{db6}(t)$ is easily derived by using the convolution theorem [7], and since that $\Psi_{long\,db6}(\omega)$ computation involves only impulses, it corresponds to a superposition of *sinc* functions (Fig. 2).

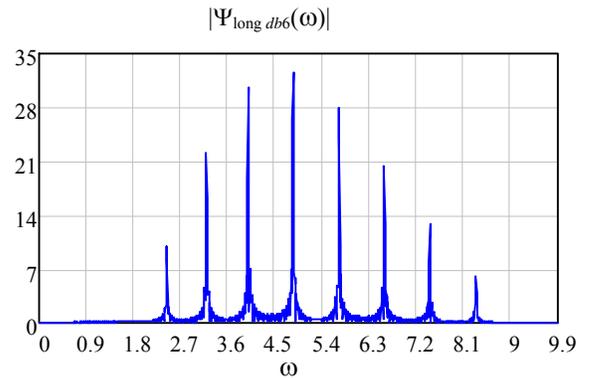

Figure 1. Spectrum of the approximated close expression, long db6, used as a model for db6. $\Psi_{longdb6}(\omega)$ is the Fourier transform pair with the inharmonic series $\psi_{long}(t)$ for db6.

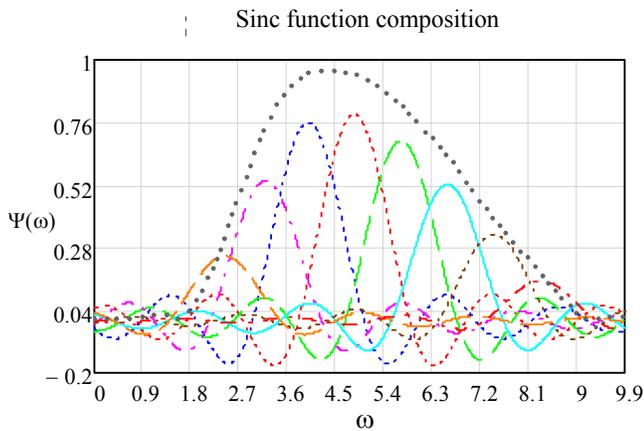

Figure 2. Spectrum of short *db*6 wave: the *sinc* pulses are superposed so as to construct the continuous wavelet spectrum, which refers to the support confinement (truncation in the time domain).

An analytical (tight) expression for the *db*6 spectrum is given by

$$\Psi_{short}(\omega) = \frac{T}{2}\sum_{k=1}^{\infty} a_k \text{sinc}\left[(\omega - b_k)\frac{T}{2\pi}\right] = |\Psi_{db}(\omega)|. \quad (16)$$

Here, T is the length of the db-wavelet support and the constants $a_k$ and $b_k$ are obtained from the approximations as described in the Table I. This is a native analytical expression for deriving the daublet spectrum. The shape of *db*6 spectrum is shown in Fig. 3.

### C. Build daublet continuous approximations in Matlab[TM] wavelet toolbox.

With the aim of investigating some potential applications of such wavelets, software to compute them according with the derived analytical close approximations should be written. Nowadays, one of the most powerful software supporting wavelet analysis is the Matlab[TM], especially when equipped with the wavelet graphic interface.

In the Matlab[TM] wavelet toolbox [12], there exists five kinds of wavelets (type the command waveinfo on the prompt):

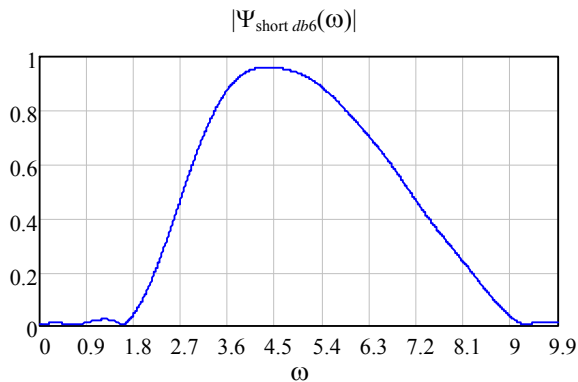

Figure 3. Result of the sinc sums, now $\Psi_{short\,db6}(\omega)$, in Fig.2: this is the same *db*6 wavelet spectrum

(i) crude wavelets; (ii) Infinitely regular wavelets; (iii) Orthogonal and compactly supported wavelets; (iv) biorthogonal and compactly supported wavelet pairs; (v) complex wavelets. Type (iii) wavelets present as a rule a scale function and possess closed expressions. However, only a few wavelets in the standard 1-D wavelet ensemble (haar, db, sym, coif, bior, rbior, meyr, dmey, gaus, mexh, morl) hold this feature. Just Meyer wavelet is assumed as type 3 wavelet, but it is not compactly supported. In spite of scaling function do exist for "dbN", these wavelets are not available as type-3 wavelets, since analytical expressions are lacking.

This paper fills this gap by providing simple analytical expressions for both wavelet and scale function of *db*N. The near-daublet expressions (for wavelet and scale functions) have been implemented in Matlab[TM] wavelet toolbox [12] as illustrated in Fig. 4 - 5, with the family name *cdb*N. In fact, *cdb*N is not a new family of wavelets, since it corresponds to *db*N. However, for the sake of convenience, *cdb*N is used as a notation to sign that continuous approximations are used in the implementation of *db*N. Potential applications of such *db*N approximation are presented in the section V.

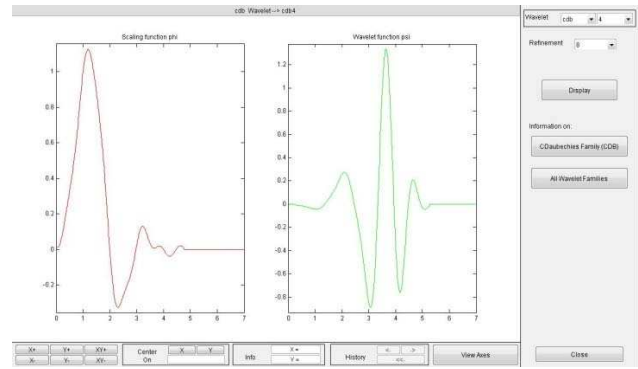

Figure 4. Example of *cdb*4 near-daublet implemented in Matlab[TM] wavelet toolbox.

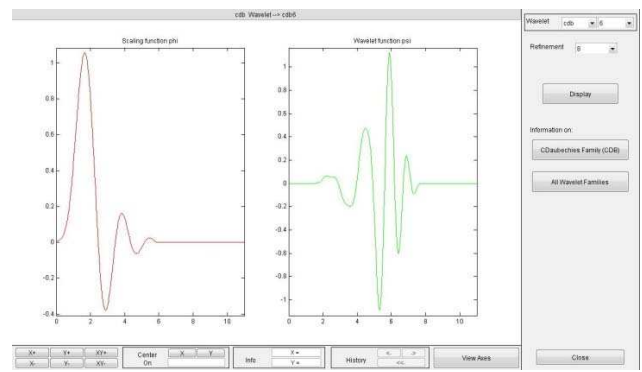

Figure 5. Display of *cdb*6 near-daublet in Matlab[TM] wavelet toolbox.

The wavelet display supplies indistinguishable waveforms as compared with those generate by the cascade algorithm (with more than 8 iterations). Here, these plots using the close expressions are inserted merely to illustrate that the close formulas derived in this paper actually work.

## IV. HAVE dbN POOR OR GOOD FREQUENCY DETECTION PROPERTIES?

The tools well suited for a spectral analysis of signals are those based on the FFT. Spectrograms are a powerful technique to investigate the frequency content in the time-frequency plane. Although wavelets are not specifically designed for carrying out spectral analysis, it can also recover some spectral information. The corresponding analysis in the wavelet scenario is performed using the scalogram, which provides location of frequency information, with wavelet-dependent accuracy [1]. On Matlab$^{TM}$, the frequency information, related to the scale, is normally retrieved by scal2frq command that computes a correspondence table of scales and frequencies (factor table scale) [12]. Some wavelets are able to detect the frequency location very well and they typically have analytical representation in the time domain. Since this paper has offered an analytical expression for *db*N, we decide to investigate the frequency detection properties of such wavelets using scalograms.

As an example, let us considerer the frequency location detection properties of *db*4 applied to a sum of two periodic signals with frequencies normalized 10 and 40. Fig. 6 plots the continuous time *db*4 wavelet analysis for that signal.

At a first glance, the evaluation provided by Fig. 6 suggested that *db*4 grants poor frequency detection properties, such as Haar (*db*2) [12]. Nevertheless, this is misleading and deserves a deeper investigation. The scalogram can be plotted in a 3D-representation as shown in Fig. 7.

In Fig. 7 the maxima energy connecting points can be examined in other viewpoint. It is observed that the 3D signal derived from the *db*4 analysis has displaced peaks in time.

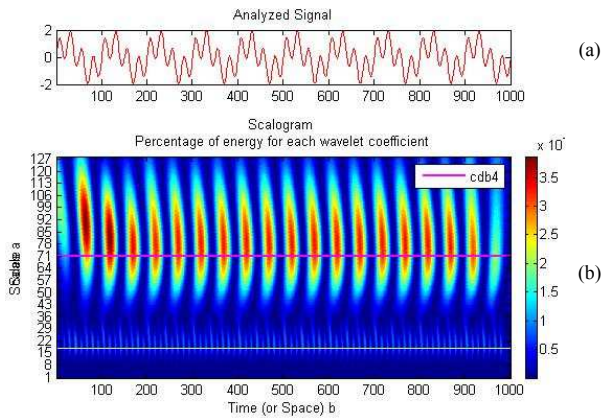

Figure 6. Plot of the connecting of maxima energy in the scalogram derived from the analysis with *cdb*4. a) analyzed signal, b) corresponding scalogram.

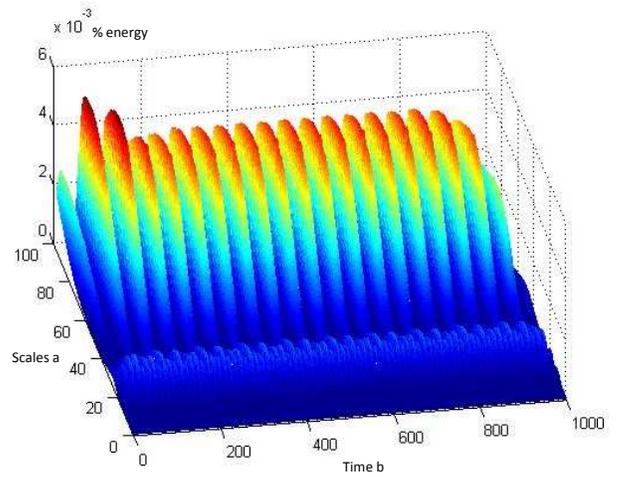

Figure 7. 3D-scalogram plot for the analysis of the signal Fig.6a with *cdb*4.

Such maxima have also a small scale displacement and the wavelet analysis ponders the lower peaks. Surprising, the most weighted group of peaks was correctly done, because they have greater amplitude. Nevertheless, due to a parallaxes error, this cannot be correctly observed in a standard 2D-scalogram such as that one shown in Fig. 6. Thus, we bring attention to the fact that 2D-scalograms seem suggest a poor detection feature for *db*N, but this is fallacious.

## V. ON POTENTIAL APPLICATIONS OF THE ANALYTICAL DAUBLET APPROACH

Despite the fact of analytical expressions for the time and spectra of continuous wavelets are known [1, 3] for both infinite supported wavelets (Shannon, Morlet, Mexhat, Meyer, de Oliveira) and compactly supported wavelets (Haar, hat, beta), very useful and powerful wavelets like the daublets (coifflets, symlets) do not present this feature. Wavelet has long been used as a powerful and deep rooted tool in signal analysis with applications in several different fields such medicine, power control, image compression, voice encoder. An extensive list of applications in Telecommunications can be found in [3]. Among them it can be mentioned a spread-spectrum system based on finite field wavelets [13], a digital wavelet-based modulation (wavelet-shift keying) [14], a multiresolution division multiplex [15]. The (AM) modulation theorem is one of the most celebrated and widely applied results of the Communication Theory [11]. By the same token as standard analog modulations, two kinds of "wavelet modulation" between a signal $f(t)$ and a continuous wavelet $\psi_{a,b}(t)$ used as a (short pulses) carrier [16] can be devised using the analytical expressions proposed in this paper.

Recent research focused on the multi-carrier transmission techniques [17, 18], highlighted that some of disadvantages inherent in OFDM systems can be steadily counteracted using wavelet carriers instead of OFDM complex exponential waveforms. Due to the fact that these wavelet carriers provide orthogonality between subcarriers and spectral containment, besides they can be separated at the receiver side by correlation techniques [19, 20]. The applications above reported, do already use wavelet carrier like Haar and Daubechies families, but need iterating process to implement that wavelet. Same processes are also required in communication channel applications such as jamming, intersymbol interference (ISI), adjacent channel interference [21].

## VI. CONCLUDING REMARKS

This paper offers a new reading of compactly supported wavelets derived from FIR filter, which was implemented through a robust algorithm and guarantee close approximations. Daublets properties of both the scaling function and the mother wavelet were maintained. Therefore, this new focus on the db wavelet family can be used to a number of applications where analytical expressions are indispensable. Time (13) and spectrum (16) close expressions for *db*N wavelets are supplied. These wavelets have been implemented on Matlab$^{TM}$ (wavelet toolbox) and an application to 3D-scalograms for detection the maxima energy peaks is shown. Arguments that *db*N provides good localization properties are presented.

This new wavelet approach seems to be particularly suitable as natural candidate to replace the daublets where they would have good performance, involving continuous wavelet-based systems, such as wavelet-OFDM, without additional computational resources to build them. As a future work, all this procedure can be naturally developed for coiflets and symlets or further discrete wavelet functions, just searching the suitable coefficients to provide the continuous approximation analytical formula.

## REFERENCES

[1]  S. Mallat, A Wavelet Tour of Signal Processing, Courant Institute, 2$^{nd}$ Ed., New York University, 1999.

[2]  A. Boggess, F.J. Narcowich. A First Course in Wavelets with Fourier Analysis, Prentice-Hall, NJ, 2002.

[3]  H.M. de Oliveira, Análise de Sinais para Engenheiros: Uma Abordagem via Wavelets, (Signal Analysis for Engineers) Editora Brasport, Rio de Janeiro, 2007.

[4]  S. Mallat, W.L. Hwang, "Singularity Detection and Processing with Wavelets," IEEE Trans. on Inform Theory, vol.38, N.2, pp.617-643, 1992.

[5]  J. Slavič, I. Simonoski, M. Boltežar, "Damping Identification using Continuous Wavelet Transform: Application to Real Data," Journal of Sound and Vibration, vol.262, N.2, pp. 291-307, 2003.

[6]  Z. Chik, T. Islam, S.A. Rosyidi, H. Sanusi, M.R. Taha, M.M. Mustafa, "Comparing the Performance of Fourier Decomposition and Wavelet Decomposition for Seismic Signal Analysis," European Journal of Scientific Research, vol.32 N.3 pp.314-328, 2009.

[7]  M. C.S. Burrus, R.A. Gopinath, H. Guo, Introduction to Wavelets and the Wavelet Transform. Prentice-Hall, Englewood Cliffs, NJ, 1998.

[8]  L. Debinath, Wavelet and Signal Processing – Applied and Numerical Harmonic Analysis, Birkhauser Boston, Berlin, 2002.

[9]  K. Levenberg, "A Method for the Solution of Certain Problems in Least Squares," Quart. Appl. Math, Vol. 2, pp. 164-168, 1944.

[10] D. Marquardt, "An Algorithm for Least Squares Estimation of Nonlinear Parameters," SIAM J. Appl. Math., Vol. 11, pp. 431-441, 1963.

[11] B.P. Lathi, and Z. Ding, Modern Digital and Analog Communication Systems, 4th Ed, Oxford University press, 2008.

[12] M. Misiti,Y. Misiti, G. Oppenheim, J.M. Poggi, Wavelet Toolbox$^{TM}$ 4 User's Guide,The MathWork Inc., MA, 2008.

[13] H.M. de Oliveira, T.H. Falk, R.G.F. Távora, "Decomposição Wavelet sobre Corpos Finitos (Wavelet decomposition over finite fields)," Rev. da Soc. Bras. de Telecomunicações, Campinas, SP, vol. 17, pp. 38 – 47, 2002.

[14] H.M. de Oliveira, H.A.N. Silva, E.A. Bouton, "Wavelet Shift-Keying: A New Digital Modulation," Proc. of the XX Simpósio Bras. de Telecomunicações, Rio de Janeiro, 5-8 Outubro, 2003.

[15] H.M. de Oliveira, and E.A. Bouton, "Multiresolution Division Multiplex (MRDM): A New Wavelet-based Multiplex System," VI Int. Telecomm. Symp. (ITS2006), September 3-6, Fortaleza, Brazil. 2006.

[16] G.A.A. Araújo, H.M. de Oliveira, "A Wavelet Modulation Theorem for Bandlimited Signals," XXVII Simpósio Brasileiro de Telecomunicações - September, Blumenau, Brazil, 2009.

[17] F. Zhao, H. Zhang, D. Yuan, "Performance of COFDM with different orthogonal Basis on AWGN and frequency Selective Channel," in Proc. of IEEE International Symposium on Emerging Technologies, Shanghai, China, pp. 473-475, 2004.

[18] Rainmaker Technologies Inc., "RM Wavelet Based PHY Proposal for 802.16.3", available on-line @ (access: 15/05/10). http://www.ieee802.org/16/tg3/contrib/802163c-01_12.pdf

[19] M. Oltean, "Wavelet OFDM Performance in Flat Fading Channels," Transactions on Electronics and Communications,vol.2, pp.52-66, 2007.

[20] D. Gupta et al., "Performance Analysis of DFT-OFDM, DCT-OFDM, and DWT-OFDM Systems in AWGN Channel," Fourth Int. Conf. on Wireless and Mobile Communications, Athens, Greece, 2008.

[21] S.W. Lee and B.H. Nam, "A Review of Wavelets for Digital Wireless Communication," Wireless Personal Communications, vol. 37, pp. 387–420, 2006.